\input amstex
\documentstyle{amsppt}
\vsize=8.7in \hsize=5.6in
\def\qed{$\square$}
\nologo
\parindent 16pt
 at 10truept \vsize=7.125in \hsize=4.5in

\leftheadtext{}
\rightheadtext{}

\topmatter
\title Extending bicolorings for Steiner triple systems
\endtitle

\author
M. Gionfriddo$^{ \star}$,  E. Guardo$^{ \star}$,  L. Milazzo$^{
\star}$
\endauthor

\dedicatory
In memory of Lucia Gionfriddo
\enddedicatory

\abstract We initiate the study of extended bicolorings of Steiner
triple systems (STS) which start with a $k$-bicoloring of an
STS($v$) and end up with a $k$-bicoloring of an STS($2v+1$)
obtained by a doubling construction, using only the original
colors used in coloring the subsystem STS($v$). By producing many
such extended bicolorings, we obtain several infinite classes of
orders for which there exist STSs with different lower and upper
chromatic number.
\endabstract

\thanks{Version: August 29, 2013}
\endthanks

\thanks {\it Mathematics Subject Classification} [2010]{Primary 05B05; Secondary 05C15, 51E10}
\endthanks

\endthanks

\thanks{$^{ \star}${Dipartimento di Matematica e Informatica, Viale A.
Doria 6, Catania - Italy}}
\endthanks

\endtopmatter

\document

\head 1.Introduction \endhead

A {\it Steiner triple system} (STS) is a pair $(V,\Cal B)$ where
$V$ is a $v$-set and $\Cal B$ is a collection of $3$-subsets of
$V$ called {\it triples} such that every $2$-subset of $V$ is
contained in exactly one triple, see \cite{4}. A {\it coloring} of
an STS $(V,\Cal B)$ is a mapping $\phi: V \rightarrow C$; the
elements of $C$ are called colors. If $\vert C \vert = k$, we have
a $k$-coloring.  For each $c \in C$, the set $\phi^{-1}(c) = \{x:
\phi(x) = c\}$ is a {\it color class}. A coloring $\phi$ of
$(V,\Cal B)$ is a {\it bicoloring} if $\vert \phi(B) \vert = 2$
for all $B \in \Cal B$. Here $\phi(B) = \bigcup_{x \in B}
\phi(x)$. Thus in  a bicoloring of $(V,\Cal B)$, every triple has
two elements in one color class and one in another class, so there
are no monochromatic triples nor polychromatic triples (i.e.
triples receiving three colors). A {\it strict} $k$-bicoloring is
one in which exactly $k$ colors are used. From now on we assume
that all our bicolorings are strict, unless the contrary is
explicitly stated.

Considerations of bicolorings of Steiner triple systems arose from
the theory of mixed hypergraphs pioneered by Voloshin
\cite{22,23}. In a mixed hypergraph setting, there are two kinds
of edges: $C$-edges which must contain two vertices colored with
the {\it same} color, and $D$-edges which must contain two
vertices of {\it different} colors. Requiring all edges of a
Steiner system to be both, $C$-triples and $D$-triples leads to
the concept of bicolorings. In the literature, often the terms
BSTS, BSQS, or bi-STS coloring are used instead of bicoloring (cf.
\cite{6,14,15,16,17,18,19}). We can also find results related to
particular color patterns for different designs in
\cite{1,5,7,9,10,11,12,13,20}.

The minimum (maximum) possible number $k$ in a strict
$k$-bicoloring of an STS is called the {\it lower} ({\it upper})
chromatic number of the STS. However, not every STS has a
bicoloring. The smallest such example occurs for STSs of order
$15$: of the $80$ nonisomorphic systems, $57$ are uncolorable. In
fact, every STS($v$) whose independence number is at most ${v
\over 3}$ is uncolorable. It is likely that almost all STSs have
this property although to best of our knowledge this remains
unproved (cf. \cite{4}).

Given a $k$-bicoloring $C$, if the cardinalities of the color classes are $n_1,n_2,\dots,n_k$,
we will write for brevity $C = C(n_1,n_2,\dots,n_k)$, and
assume, unless stated to the contrary, that $n_1 \leq n_2 \leq \dots \leq n_k$.

\vskip 10pt

In this paper, we want to initiate a study of {\it extended} bicolorings, i.e. bicolorings of an STS($w$) which start
with a bicoloring of a sub-STS($v$).
Essential for us in this endeavor will be a well-known recursive construction known as a {\it doubling construction}
(other names: $v \rightarrow 2v+1$ rule,
doubling plus one construction etc.) which starts with an STS($v$) and ends with an STS($2v+1$).

To obtain such a construction, all that is needed, apart from the
subsystem, is a $1$-factorization of the complete graph $K_{v+1}$.
Indeed, let $(X,\Cal F)$ where $\Cal F =
\{F_1,\dots,F_v\}$ is a $1$-factorization of $K_{v+1}$ (where
$\vert X \vert = v+1$ must be even). If $(V,\Cal B)$, $V =
\{a_1,\dots,a_v\}$, is an STS($v$), form the set of triples $\Cal
C = \{\{a_i,x,y\}: a_i \in V, \{x,y\} \in F_i\}$. Then $(V \cup
X,\Cal B \cup \Cal C)$ is an STS($2v+1$) (cf. \cite{4}).

An easy observation is that if a given STS($v$), $(V,\Cal B)$,
admits a $k$-bicoloring $C = C(n_1,\dots,n_k)$, then {\it any}
STS($2v+1$) obtained from $(V,\Cal B)$ by a doubling construction
admits a $(k+1)$-bicoloring $C(n_1,\dots,n_k,n_{k+1})$ where the
$v+1$ vertices of $X$ are colored with a new color, and so $
n_{k+1} = v+1$. Another such $(k+1)$-bicoloring that can be always
obtained is $C' = C'(n_1',\dots,n_k',n_{k+1}')$ where $n_i' =
2n_i$ for $i=1,\dots,k$, and $n_{k+1}' = 1$ (see \cite{3}).

\vskip 10pt

The question that we want to address is the following. Given an
STS($v$) with a bicoloring $C = C(n_1,\dots,n_k)$, when does
there exist an STS($2v+1$) obtained by a doubling construction
which admits a bicoloring $C'  = C'(n_1',\dots,n_k')$? In other
words, when can we color the elements of $X$ with the original $k$
colors of the $k$-bicoloring $C$ {\it without} introducing an
extra color as above?  If such a coloring exists, we call it an
{\it extended} bicoloring of $C$. Thus extended bicolorings may
exist only for orders $2v+1 \equiv 3$ or $7\ (mod\ 12)$ as $v
\equiv 1$ or $3\ (mod\ 6)$.

\vskip 5pt

The importance of extended bicolorings lies in the fact that they
enable one to construct STSs with different lower and upper
chromatic numbers; there are only scarce results in the literature
on the latter (cf., e.g., \cite{14}). The extendibility of partial
colorings is a relevant issue in graph theory (see \cite{2,21}),
both for its theoretical interest and practical applications. Here
we initiate the study of extending bicolorings in Steiner triple
systems, with the aim to derive consequences in the coloring
theory of mixed hypergraphs. In this way the present work relates
several intensively studied areas.

\head 2. Extended bicolorings \endhead

Let $S = (V,\Cal B)$ be an STS($v$) which is $k$-bicolorable with
$C = C(n_1,\dots,$  \newline $n_k)$ and let $S' = (X,\Cal C)$ be
an STS($2v+1$) obtained  from $S$ by a doubling construction.  We
are trying to investigate the conditions under which there exists
an {\it extended} bicoloring of $S'$, say $C' =
C'(n_1',\dots,n_k')$ where the elements of the subsystem $(V,\Cal
B)$ are colored as in $C$, and the elements of $Y =X \setminus V$
are colored with the same colors as those used in $C$. If $c_i =
n_i' -n_i$, $1 \leq i \leq k$ are the numbers of vertices in  $Y$
colored with the color $i \in C$ then clearly, $\sum_{i=1}^{k} c_i
= v+1$ (it may happen that $c_j = 0$ for some $j \in
\{1,\dots,k\})$. Beside this obvious condition, the following is a
necessary condition for the existence of an extended
$k$-bicoloring of $S'$.

\vskip 10pt

\proclaim{Theorem 1} Let $S = (V,\Cal B)$ be an STS($v$) which is
$k$-bicolorable with $C = C(n_1,\dots,n_k)$ and let $S' = (X,\Cal
C)$ be an STS($2v+1$) obtained  from $S$ by a doubling
construction. With the notation as above,

\vskip 10pt
\hskip 50pt
$\sum_{i=1}^{k} c_i^2 + 2 \sum_{i=1}^{k}n_ic_i = (v+1)^2$.    \hskip 100pt  $(1)$
\endproclaim

\noindent {\bf Proof.} The number of pairs of elements of $Y$ equals ${{v+1} \choose 2}$.
Clearly, the number of monochromatic pairs among these is
$\sum_{i=1}^{k} {c_i \choose 2}$. On the other hand,
the number of two-colored pairs among these is
$\sum_{i=1}^{k} n_ic_i$; indeed, if $a_i \in V$ is colored with color $j$, then any pair $\{x,y\}$ in the $1$-factor $F_i$ is either monochromatic or else one of $x,y$ is colored with color $j$. Consequently, the number of two-colored pairs in $F_i$ is $n_j$. Thus we have

\vskip 5pt

\hskip 20pt  $\sum_{i=1}^{k} {c_i \choose 2} \ + \ \sum_{i=1}^{k} n_ic_i = {{v+1} \choose 2}$

\vskip 10pt

\noindent whence (1) follows easily.                      \qed

\vskip 10pt

Solutions of (1) will be called {\it solutions with respect to
$C$}. We stress that condition (1) given by Theorem 1 is only
necessary for the existence of an extended bicoloring. It
certainly  is not sufficient: in \cite{6} condition (1) was
determined for $v=2^h-1$ and all of its solutions  were determined
for $h \leq 10$, nevertheless these solutions do not lead to any
extended bicolorings.

\vskip 10pt

\proclaim{Corollary 2}  Let $S'$ be a $k$-bicolorable STS($2v+1$)
obtained by a doubling construction from a $k$-bicolorable
STS($v$) with the coloring $C = C(n_1,\dots,$ \newline $n_k)$, and
let $(c_1,\dots,c_k)$ be a solution to (1) with respect to $C$. \newline
1. If  there is a $c_j = 0,$ then all $c_i$'s are even. \newline
2. If there is $j$ such that $c_j > {{v+1} \over 2},$ then there
exists no extended bicoloring of $C$.
\endproclaim
\vskip 10pt

{\bf Proof.} 1. If there is a $j$ such that $c_j=0,$ then in any
factor corresponding to an element $a_l \in V$ colored with the
color $j$, all pairs must be monochromatic which implies that
every $c_i$ has to be even. \newline 2. If $c_j > {{v+1} \over 2}$
for some $j,$ then in all factors associated with elements of $V$
colored with color $j$, there must exist monochromatic pairs of
color $j$, and thus monochromatic triples, which is a
contradiction.     \qed

\vskip 10pt

We illustrate the use of Corollary 2 on the example of a
(potential) extended bicoloring of an STS($19$). First notice that
no extended bicolorings of STS($v$) exist for $v=7$ or $v=15$, as
shown in \cite{6}. The unique STS($9$) admits a bicoloring $C =
C(1,4,4)$, and no other bicolorings (see \cite{3} or \cite{17}).
The following are all solutions with respect to $C$: (a)
$(3,2,5)$, (b) $(3,5,2)$, (c) $(5,0,5)$, (d) $(5,5,0)$, (e)
$(8,0,2)$, and (f) $(8,2,0)$. Corollary 2.1 eliminates solutions
(c), (d), (e) and (f) from contention. Concerning (a), since
$c_1=3$ and $c_3=5$, it must be that in the four $1 $-factors
associated with elements colored with color $2$, there are exactly
two $2$-colored pairs colored with colors $1$ and $2$ and with $2$
and $3$, one monochromatic pair of color $1$, and two
monochromatic pairs of color $3$. Since $c_1=3$, this is easily
seen to be impossible, so solution (a) cannot lead to an extended
bicoloring. The same reasoning applies to the solution (b). Thus
there exist no extended bicolorings of any STS($19$) (obtained
from an STS($9$) by a doubling construction, of course). Thus the
smallest $w$ for which an  STS($w$) may admit an extended
bicoloring is $w=27$ (where the STS($27$) is obtained from an
STS($13$) by a doubling construction).

We remark that due to the above, any uniquely $3$-colorable
STS($19$), or any $3$- and $4$-colorable STS($19$) cannot contain
a sub-STS($9$). It was shown in \cite{14} that there exist
uniquely $3$-bicolorable, uniquely $4$-bicolorable, and also $3$-
and $4$-bicolorable STS($19$).

\proclaim{Theorem 3} Let $S$ be a $k$-bicolorable STS($v$) with the $k$-bicoloring $C = C(n_1,\dots,n_k)$. and suppose there
exist $i, j, i \neq j$, such that
$n_i + n_j = {{v+1} \over 2} \equiv 0\ (mod\ 2)$. Then there exists an  STS($2v+1)$, $S'$, obtained by a doubling construction from $S$ such
that $S'$ has an extended $k$-bicoloring $C'$.
\endproclaim

\noindent {\bf Proof.} Since $v+1 \equiv 0\ (mod\ 4)$, we may use
as the $1$-factorization $\Cal F = \{F_1,\dots,F_v\}$ in the
doubling construction the following $1$-factorization. Write $Y =
Y_1 \cup Y_2$ where $\vert Y_i \vert = {{v+1} \over 2}$; take
$F_1,\dots,F_{{v+1} \over 2}$ to be the $1$-factors of any
$1$-factorization of the complete bipartite graph $K_{{{v+1} \over
2},{{v+1} \over 2}}$ with bipartition $(Y_1,Y_2)$; for the
remaining ${{v-1} \over 2}$ $1$-factors $F_{{v+3} \over
2},\dots,F_v$, take  $F_i = G_i \cup H_i$, $i={{v+3} \over
2},\dots,v$, where $G_i, H_i$ are the $1$-factors of any
$1$-factorization of $K_{{v+1} \over 2}$ on $Y_1$, and $Y_2$,
respectively. Color now the ${{v+1} \over 2}$ vertices of $Y_1$
with color $i$ and the ${{v+1} \over 2}$ vertices of $Y_2$ with
color $j$. Associate the $1$-factors $F_1,\dots,F_{{v+1} \over 2}$
with the vertices of $V$ colored in the coloring $C$ with either
color $i$ or color $j$, and associate the remaining $1$-factors
with the elements of $V$ colored in $C$ with colors other than $i$
or $j$. We obtain in this way an extended $k$-bicoloring of the
resulting STS($2v+1)$. Indeed, if $a_q$ is an element of $V$ which
is colored with $i$ or $j$, then any triple $T$ containing $a_q$
is two-colored: one of the two elements of $T$ other than $a_q$ is
colored with color $i$, and the other with color $j$. On the other
hand, if $a_r$ is an element of $V$ colored in $C$ with a color
other than $i$ or $j$, then any triple $T$ containing $a_r$ is
also two-colored since the two elements of $T$ other than $a_r$
are {\it both} colored with $i$ or both colored with $j$. \qed

A more general version of Theorem 3 is the following.

\proclaim{Theorem 4} Let $S$ be a $k$-bicolorable STS($v$) with the $k$-bicoloring $C = C(n_1,\dots,n_k)$. Suppose that there exist $p$ integers
$n_{k_i}$, $1 \leq i \leq p < k$ such that $n_{k_1} + n_{k_2} = {{v+1} \over {2^{p-1}}}$ is an even integer, and further $n_{k_i} = {{v+1} \over {2^{p-i+1}}}$ for $3 \leq i \leq p$ are all even. Then there exists an STS($2v+1$) obtained by a doubling construction from $S$ which
has an extended $k$-bicoloring.
\endproclaim

The proof of this theorem is more technical than that of Theorem
3, especially in the description of the $1$-factorization $\Cal F$
involved in the doubling construction. Since in what follows we do
not make use of this more general version, with one exception,
this proof is omitted (see Appendix \cite{8}).

\head 3. Small extended bicolorings \endhead

As shown earlier, there exist no extended bicolorings of STS($w$)
for $w \leq 19$. Since $w \equiv 3$ or $7\ (mod\ 12)$, the
smallest $w$ for which there might exist an extended bicoloring is
$w=27$. Such an extended bicoloring does indeed exist.

\proclaim{Theorem 5} There exists an STS($27$), $(W,\Cal C)$
obtained by a doubling construction from an STS($13$), $(V,\Cal
B)$, which has an extended $3$-bicoloring  $C=C(2,5,6)$. For this
system, $\chi =3$ and $\bar{\chi} = 4$.
\endproclaim

\noindent {\bf Proof.} All solutions $(c_1,c_2,c_3)$ with respect
to the coloring $C$ (cf. Theorem 1) are as follows: (a) $(4,4,6)$,
(b) $(7,1,6)$, (c) $(4,7,3)$, (d) $(7,7,0)$, (e) $(10,1,3)$, (f)
$(10,4,0)$.  By Corollary 2, solutions (d), (e), and (f) cannot
lead to an extended bicoloring of $C$. Concerning solution (c),
there are three monochromatic pairs of elements of color $3$. Two
of these pairs may occur in the two $1$-factors corresponding to
the vertices of $V$ of color $1$ but the third pair cannot occur in
a $1$-factor corresponding to a vertex of $V$ of color $2$ (as
there are $7$ vertices of $W \setminus V$ of color $2$, and that
would force a monochromatic triple of color $2$), nor clearly in a
$1$-factor corresponding to a vertex of color $3$. Thus solution
(c) does not lead to an extended bicoloring of $C$ either.

On the other hand, each of the first two solutions, namely
$(4,4,6)$ and $(7,1,6)$, lead to an extended bicoloring
$C'=C'(6,9,12)$. The $1$-factorizations $\Cal F$ used in the
respective doubling constructions are given in the Appendix
\cite{8}.  Our STS($27$), $(W,\Cal C)$,  besides having an
extended $3$-bicoloring with respect to $C$, is also
$4$-bicolorable with the coloring $C"=C"(2,5,6,14)$. At the same
time, a $5$-bicoloring of $(W,\Cal C)$  is impossible due to
\cite{20}, since $27 < 2^5-1$. Thus $\chi =3, \bar{\chi} = 4$,
as claimed. \qed

\vskip 10pt

Concerning order $31$, an inspection of the tables in \cite{3}
shows that there exists no extended bicoloring for this order:
there exists no $3$-bicoloring of an STS($15$), and no
$4$-bicoloring of STS($31)$ whatsoever. However, the next
admissible order  $39$ shows a quite different behaviour.

\proclaim{Theorem 6} There exist STS($39$) admitting extended
bicolorings obtained from extended bicolorings of STS($19$) of
type $C_1=C(4,6,9)$ and $C_2=C(1,2,8,8)$. More specifically, there
exist STS($39$) with $(\chi,\bar{\chi})$ equal to either $(3,4)$,
or $(4,5)$, or $(3,5)$.
\endproclaim

\noindent {\bf Proof.} It was shown in \cite{14} that there exist
STS($19$) (a)  admitting only the $3$-bicoloring $C(4,6,9)$, (b)
admitting only the $4$-bicoloring $C(1,2,8,8)$, and (c) admitting
both, the $3$-bicoloring $C(4,6,9)$ and the $4$-bicoloring
$C(1,2,8,8)$.  By Theorem 3, both $C_1$ and $C_2$ are extendable
bicolorings (since we have $4+6 = 10$ and $2+8 = 10$,
respectively). Starting with an STS($19$) of type (a), (b), or
(c), we obtain an STS($39$) of the respective kind as claimed.
\hskip 30pt        \qed

 \vskip 10pt

In what follows we discuss in  somewhat less detailed manner the existence of extended bicolorings for those STS($v$) of orders $43\leq v \leq 99$
which can  be obtained by a doubling construction.

\proclaim{Theorem 7} For an STS($v$), $v \in \{51,63,67,75\}$, there exists no extended bicoloring.
\endproclaim

\noindent {\bf Proof.} (i) For an STS($25$), the only types of
bicoloring that are possible are $C_1 = C(5,10,10)$ and $C_2 =
C(1,4,8,12)$. While there exist $12$ solutions with respect to
$C_1$, none satisfies the condition of Corollary 2 and thus cannot
lead to an extended bicoloring. There exist no solutions with
respect to $C_2$, and so no STS($51$) can have an extended
bicoloring. \newline (ii) By \cite{6}, no STS($63$) obtained by
doubling from  an STS($31$) can have an extended bicoloring.
\newline (iii) None of the bicolorings of any STS($33$) or
STS($37$) (cf. \cite{3}) yields a solution with respect to such a
coloring, thus there is no extended bicoloring of any STS($67$) or
STS($75$).   \hskip 20pt  \qed

\vskip 10pt

\proclaim{Theorem 8} For an STS($v$), $v \in
\{43,55,79,87,91,99\}$, there exist extended bicolorings. More
specifically,  there exists an extended $3$-bicoloring of an
STS($43$), extended $4$-bicolorings of an STS($v$) for $v \in
\{55, 87, 91\}$, and extended $4$- and $5$-bicolorings of an
STS($v$) for $v \in \{79, 99\}$.
\endproclaim

\noindent {\bf Proof.} (i) A bicolorable STS($21$) can only be
$3$-bicolorable, with colorings $C_1=C(5,6,10)$ or $C_2=C(4,8,9)$
(see \cite{3} or \cite{14}). Both are extendable to a
$3$-bicoloring $C = C(10,16,17)$ of an STS($43$) for which we have
$\chi = 3$ and $\bar{\chi}=4$ (there exists no $5$-bicolorable
STS($43$), cf. \cite{3}). The $1$-factorization $\Cal F$ in the
corresponding doubling construction is given in the Appendix
\cite{8}. \newline (ii) By Theorem 3, the $4$-bicoloring
$C(1,\bold{4},\bold{10},12)$ of an STS($27$) is extendable to a
$4$-bicoloring $C(1,12,18,24)$ of an STS($55$), further the
$4$-bicoloring $C(1,\bold{8},\bold{12},18)$, and the
$4$-bicoloring $C(\bold{2},6,13,\bold{18})$, respectively, of an
STS($39$) is extendible to a $4$-bicoloring $C(1,18,28,32)$, and
to a $4$-bicoloring $C(6,13,22,38)$ of an STS($79$), respectively;
finally, the  $4$-bicolo-ring $C(1,\bold{10},\bold{12},20)$, and
the $4$-bicoloring $C(\bold{4},4,17,\bold{18})$, respectively, of
an STS($43$) is extendable to a $4$-bicoloring $C(1,20,32,34)$,
and to a $4$-bicoloring $C(4,17,26,40)$, respectively. (The two
essential colors are indicated in bold.) \newline (iii)
Extendability of the $5$-bicolorings $C(1,2,8,8,20)$ and
$C(1,4,4,10,20)$ of an STS($43$) follows from the more general
Theorem 4. \newline (iv) There are only two possible types of a
$4$-bicoloring of an STS($45$), namely $C_1=C(2,8,14,21)$ and
$C_2=C(4,6,13,22)$. There are $12$ solutions with respect to $C_1$
but none of them leads to an extended bicoloring. Similarly, there
are $12$ solutions with respect to $C_2$ but only one of them,
namely $(c_1,c_2,c_3,c_4) = (4,8,12,22)$ leads to an extended
$4$-bicoloring. The corresponding $1$-factorization $\Cal F$ in
the doubling construction that leads to this extended bicoloring
is given in the Appendix \cite{8}. \newline (v) Although there
exist $3$-, $4$-, and $5$-bicolorable STS($49$), none of the
$3$-bicolorings is extendable. On the other hand, $4$-bicolorings
$C(2,8,18,21)$ and $C(5,6,14,24)$ as well as the $5$-bicoloring
$C(1,4,4,20,20)$ are all extendable. This is shown by examining
all solutions with respect to the particular bicoloring $C$. Due
to the considerable number of these solutions ($84$, $29$ and
$27$, respectively), we omit the details. The $1$-factorizations
occurring in the doubling constructions leading to the respective
extended bicolorings are given in the Appendix \cite{8}.

\vskip 10pt

\proclaim{Theorem 9} There exist extended $4$-bicolorings for each order $w \in \{127,$ \newline $151,159,175\}$; there exist extended $5$-bicolorings for
each order $w \in \{103,$ \newline $111,127,135,151,159,175\}$.
\endproclaim

\noindent {\bf Proof.} Below we list $4$- and $5$-bicolorings
(known to exist by \cite{3}) to which it is possible to apply
Theorem 3; the two essential colors are in bold.

\vskip 10pt

\noindent order $2v+1$    \hskip 20pt  extendable colorings of an STS($v$)

$103          \hskip 50pt        (1,\bold{2},8,16,\bold{24})$

$111         \hskip 50pt         (1,2,\bold{8},\bold{20},24)$

$127       \hskip 50pt          (2,\bold{14},\bold{18},29), \ (\bold{4},9,22,\bold{28}), \  (\bold{2},5,6,20,\bold{30})$

$135     \hskip 50pt          (1,\bold{2},16,16,\bold{32})$

$151    \hskip 50pt           (4,\bold{12},\bold{26},33), \  (1,4,\bold{10},\bold{28},32)$

$159   \hskip 50pt           (\bold{4},14,25,\bold{36}), \ (\bold{6},10,29,\bold{34}), \ (1,\bold{4},12,26,\bold{36}),$ \newline

\hskip 63pt      $ \ (1,2,\bold{16},\bold{24},36)$

$175  \hskip 50pt     (\bold{4},17,26,\bold{40}, \ (2,5,\bold{10},\bold{34},36)$.

\vskip 10pt

We summarize our results as follows.

\proclaim{Theorem 10} Let $\Omega = \{27,39,43,55,79,87,91,99,103,111,127,135,151,$ \newline $159,175\}$. For each $v \in \Omega$, there exists an STS($v$)
with an extended bicoloring, and thus for all $v \in \Omega$, we have $\chi \neq \bar{\chi}$.
\endproclaim

\noindent {\bf Proof.} For each $v \in \Omega$ we have an extended $k$-coloring for some $k$, and also (at least) a $(k+1)$-bicoloring (with $n_{k+1} = v+1)$.
\qed

\vskip 10pt

\proclaim{Corollary 11} For each $v \in \Omega' =
\{27,39,43,91,99,103,127,135,151\}$, there exists an infinite
class of STS($w$), where $w = 2^t(v+1)-1$, $t \geq 1$, such that
$\chi \neq \bar{\chi}$.
\endproclaim

\noindent {\bf Proof.} Apply repeatedly the doubling construction to the appropriate STS($v$).   \qed

\head 4. Conclusion \endhead

In  this paper, we have investigated extended bicolorings with the
explicit aim to prove the existence of STSs with $\chi \neq
\bar{\chi}$, that is, with different lower and upper chromatic
number. We established the existence of extended bicolorings and
of such STSs for several infinite classes of orders $2v+1 \equiv
3$ or $7\ (mod\ 12)$, by utilizing the doubling construction. The
problem of determining for which orders $v \equiv 1$ or $3\ (mod\
6)$ does there exist an STS($v$) with different lower and upper
chromatic number is certainly worthwhile. Another interesting
question is, how large can the difference $\bar{\chi} \ - \ \chi$
be? It is also a legitimate question to ask whether an analogue of
extended bicolorings may exist for other recursive constructions,
such as the known $v \rightarrow 2v+t$ rules where $t > 1$ (cf.
\cite{4}). For example, is it possible to use the $v \rightarrow
2v+5$ rule starting with an STS($7$) and ending up with an
STS($19$) to show that the $3$-bicoloring $(1,2,4)$ for STS($7$)
can be extended to a $3$-bicoloring $(4,6,9)$ for an STS($19$)?
The next example answers this question.

\vskip 5pt

\noindent {\bf Example 12.} The following STS($19$) with a
sub-STS($7$) has a $3$-bicoloring and these triples: $\{0,1,9\},
\{2,3,9\}, \{0,2,10\}, \{1,3,10\}, \{0,3,15\},$
\newline $\{1,2,15\}, \{9,10,15\}, \{0,4,11\}, \{0,5,12\},
\{0,6,13\}, \{0,7,14\}, \{0,8,16\},$
\newline $\{1,4,12\}, \{1,5,11\}, \{1,6,14\}, \{1,7,13\}, \{1,8,17\}, \{2,4,13\}, \{2,5,14\}$, \newline $\{2,6,11\}, \{2,7,12\}, \{2,8,18\}, \{3,4,14\},
\{3,5,16\}, \{3,6,17\}, \{3,7,18\}$, \newline $\{3,8,11\}, \{4,5,9\}, \{4,6,18\}, \{4,7,16\}, \{4,8,10\}, \{5,6,10\}, \{5,7,17\},$ \newline $\{5,8,13\},
\{6,7,9\}, \{6,8,12\}, \{7,8,15\}, \{9,11,16\}, \{9,12,17\}, \{9,13, 18\},$ \newline $\{8,9,14\}, \{7,10,11\}, \{10,12,16\}, \{10,13,17\}, \{10,14,18\},
\{11,12,18\},$ \newline $\{11,13,15\}, \{11,14,17\}, \{3,12,13\}, \{12,14,15\}, \{13,14,16\}, \{6,15,16\},$ \newline $\{4,15,17\}, \{5,15,18\}, \{2,16,17\},
\{1,16,18\}, \{0,17,18\}$.

\noindent The first seven triples are those of an STS($7$) on
$\{0,1,2,3,9,10,15\}$; the three color classes are
$\{0,1,2,3,4,5,6,7,8\}, \{9,10,11,12,13,14\}$ and \newline
$\{15,16,17,18\}$.

\vskip 10pt

Even if the answer in this case proved to be affirmative, and may
proved so in similar cases, it is not immediately clear that this
will have as a consequence the existence of STSs with $\chi \neq
\bar{\chi}$. Thus the doubling construction appears to offer most
benefits from the stated applications point of view. Nevertheless,
it seems to us worthwhile to study ``extended" bicolorings for
recursive rules for STSs other than doubling.

\head Acknowledgements \endhead

Thanks to Alex Rosa for valuable comments, and to Mariusz Meszka
for providing us with Example 12. The authors also thank the
reviewers and editors for their assistance.

\enddocument